\documentclass[oneside,english]{amsart}
\usepackage[T1]{fontenc}
\usepackage[latin1]{inputenc}
\usepackage{amssymb}
\IfFileExists{url.sty}{\usepackage{url}}
                      {\newcommand{\url}{\texttt}}

\makeatletter
 \theoremstyle{plain}    
 \newtheorem{thm}{Theorem}[section]
 \numberwithin{equation}{section} 
 \numberwithin{figure}{section} 
 \theoremstyle{plain}
 \theoremstyle{plain}    
 \newtheorem{prop}[thm]{Proposition} 
 \theoremstyle{plain}    
 \newtheorem{lem}[thm]{Lemma} 

\usepackage{babel}
\makeatother
\begin{document}

\title{Poisson actions up to homotopy and their quantization}

\author{Pavol \v Severa}

\begin{abstract}
Symmetries of Poisson manifolds are in general quantized just to symmetries
up to homotopy of the quantized algebra of functions. It is therefore
interesting to study symmetries up to homotopy of Poisson manifolds.
We notice that they are equivalent to Poisson principal bundles and
describe their quantization to symmetries up to homotopy of the quantized
algebras of functions, using the formality theorem of Kontsevich.
\end{abstract}

\address{Section de Mathématiques, Univ.~de Gen\`eve, 2-4 rue du Li\`evre,
1211 Gen\`eve, Switzerland, on leave from Dept.~of Theoretical Physics,
Mlynská dolina F2, 84248 Bratislava, Slovakia}

\thanks{supported in part by the Swiss National Science Foundation}

\maketitle
An action of a Lie algebra $\mathfrak{g}$ on a Poisson manifold $M$ is
\emph{Poisson}, if it preserves the Poisson structure on $M$. In
general it is not possible to quantize $M$ so that $\mathfrak{g}$ would
act by derivations on the quantized algebra of functions. On the other
hand, as a simple application of the formality theorem of Kontsevich
\cite{K1}, a Poisson action up to homotopy (a generalization of Poisson
actions, defined below) is easily quantized to an action up to homotopy
on the quantized algebra of functions. It is therefore interesting
to study the geometry of Poisson actions up to homotopy, together
with its generalizations where $\mathfrak{g}$ is replaced by a Lie bialgebra.
Namely, it turns out that an up to homotopy Poisson action is the
same as a Poisson principal bundle. This geometry sheds a new light
even on some pretty standard things, e.g.~relations between Poisson
and Hamiltonian actions.

\section{Up to homotopy Poisson actions\label{sec:HPA}}

An \emph{up to homotopy Poisson action} of a Lie algebra $\mathfrak{g}$
(\emph{$\mathfrak{g}$-HPA} for short) on a Poisson manifold $M$ is an
extension $\mathfrak{g}_{M}$ of $\mathfrak{g}$ by the Lie algebra $C^{\infty}(M)$
(with the Poisson bracket used as the Lie bracket), such that for
every $X\in\mathfrak{g}_{M}$, the map $C^{\infty}(M)\rightarrow C^{\infty}(M)$,
$f\mapsto[X,f]$, is a derivation (i.e.~a vector field on $M$).
A \emph{morphism} from a $\mathfrak{g}$-HPA on a Poisson manifold $M_{1}$
to a $\mathfrak{g}$-HPA on $M_{2}$ is a Poisson map $f:M_{1}\rightarrow M_{2}$
and a morphism of the extensions of $\mathfrak{g}$ going the opposite direction,
equal to $f^{*}:C^{\infty}(M_{2})\rightarrow C^{\infty}(M_{1})$ on the kernels.

A $\mathfrak{g}$-HPA on a point is a central extension of $\mathfrak{g}$ by
$\mathbb{R}$. Given a true Poisson action of $\mathfrak{g}$ on a Poisson
manifold $M$, the corresponding extension is the semidirect product.
Any splitting of a $\mathfrak{g}$-HPA to a semidirect product turns it
into a true Poisson action ($\mathfrak{g}_{M}$ has a Poisson action on
$M$ by definition, a splitting $\mathfrak{g}\rightarrow\mathfrak{g}_{M}$ gives us then
a Poisson action of $\mathfrak{g}$). In other words, an ordinary Poisson
action on $M$ is the same as a $\mathfrak{g}$-HPA $\mathfrak{g}_{M}$ on $M$
\emph{and} a morphism from $\mathfrak{g}_{M}$ to the trivial $\mathfrak{g}$-HPA
on a point.

Let us give an equivalent definition using the Maurer-Cartan equation.
If $M$ is a manifold, let $L_{M}$ denote the differential graded
Lie algebra (DGLA) of mutivector fields $\Gamma(\bigwedge TM)[1]$
(with Schouten bracket and zero differential). A $\mathfrak{g}$-HPA on
$M$ is a degree-1 element $\sigma$ of the DGLA $\bigwedge\mathfrak{g}^{*}\otimes L_{M}$,
satisfying the Maurer-Cartan (MC) equation\begin{equation}
d\sigma+[\sigma,\sigma]/2=0.\label{eq:MC}\end{equation}
As we'll see easily, this is indeed equivalent to the previous definition,
when we choose a splitting of the extension $\mathfrak{g}_{M}=\mathfrak{g}\oplus C^{\infty}(M)$
(direct sum of vector spaces -- not direct (or semidirect) sum of
Lie algebras!). 

Let us decompose $\sigma$ into bihomogeneous parts \begin{equation}
\sigma=\sigma^{0}+\sigma^{1}+\sigma^{2}\label{eq:dec}\end{equation}
 (the superscript is the degree in $\bigwedge\mathfrak{g}^{*}$). Then $\sigma^{0}$
is a Poisson structure on $M$, and we get a Lie bracket on $\mathfrak{g}_{M}=\mathfrak{g}\oplus C^{\infty}(M)$
via $[u,v]=[u,v]_{\mathfrak{g}}+\sigma^{2}(u,v)$, $[u,f]=\mathcal{L}_{\sigma^{1}(u)}f$,
$[f,g]=\{ f,g\}$ (for any $u,v\in\mathfrak{g}$ and $f,g\in C^{\infty}(M)$).

Yet another equivalent formulation -- we have an $L_{\infty}$-morphism
from $\mathfrak{g}$ to the DGLA $L_{M,\sigma^{0}}$ (where $L_{M,\sigma^{0}}$
is the same as $L_{M}$, but with differential equal to $\mathit{ad}_{\sigma^{0}}$).

Just for the fun of it, let us plug the decomposition (\ref{eq:dec})
into the Maurer-Cartan equation (\ref{eq:MC}) and write all the results
explicitly:

\begin{enumerate}
\item $[\sigma^{0},\sigma^{0}]=0$, i.e.~$\sigma^{0}$ is a Poisson structure
on $M$
\item $[\sigma^{0},\sigma^{1}]=0$, i.e.~$\sigma^{1}\in\mathfrak{g}^{*}\otimes\Gamma(TM)$
maps elements of $\mathfrak{g}$ to vector fields preserving the Poisson
structure $\sigma^{0}$
\item $d\sigma^{1}+[\sigma^{1},\sigma^{1}]/2+[\sigma^{0},\sigma^{2}]=0$,
i.e.~$\sigma^{1}$ is not necessarily an action of $\mathfrak{g}$ on $M$,
but rather for any $u,v\in\mathfrak{g}$ we have \[
[\sigma^{1}(u),\sigma^{1}(v)]=\sigma^{1}([u,v])+X_{\sigma^{2}(u,v)},\]
 where $X_{f}$ denotes the Hamiltonian vector field generated by
function $f$
\item $d\sigma^{2}+[\sigma^{1},\sigma^{2}]=0$, i.e.~\[
\mathcal{L}_{\sigma^{1}(u)}\sigma^{2}(v,w)+\sigma^{2}(u,[v,w])+\mathit{c.p.}=0\]
where {\it c.p.} means cyclic permutations in $u,v,w$.

\end{enumerate}
Finally, let us see what happens to $\sigma$ when we choose a different
splitting of $\mathfrak{g}_{M}$ to $\mathfrak{g}\oplus C^{\infty}(M)$. The resulting
\emph{gauge transformation} of $\sigma$ by a $\tau\in\mathfrak{g}^{*}\otimes C^{\infty}(M)$
is explicitly:

\begin{enumerate}
\item $\sigma^{0}\mapsto\sigma^{0}$, i.e.~the Poisson structure on $M$
doesn't change
\item $\sigma^{1}\mapsto\sigma^{1}+[\tau,\sigma^{0}]$, i.e.~$\sigma^{1}(u)\mapsto\sigma^{1}(u)+X_{\tau(u)}$
\item $\sigma^{2}\mapsto\sigma^{2}+d\tau+[\tau,\sigma^{1}]+[\tau,[\tau,\sigma^{0}]]/2$,
i.e.~\[
\sigma^{2}(u,v)\mapsto\sigma^{2}(u,v)+\tau([u,v])+\mathcal{L}_{\sigma^{1}(u)}\tau(v)-\mathcal{L}_{\sigma^{1}(v)}\tau(u)+\{\tau(u),\tau(v)\}/2\]

\end{enumerate}

\section{Up to homotopy actions as Poisson principal bundles}

\begin{prop}
The category of $\mathfrak{g}$-HPA's is equivalent to the category of Poisson
principal $\mathfrak{g}^{*}$-bundles. Given a Poisson principal $\mathfrak{g}^{*}$-bundle
$P\rightarrow M$, we construct the corresponding $\mathfrak{g}$-HPA on $M$ as
follows: the elements of $\mathfrak{g}_{M}$ over $u\in\mathfrak{g}$ are functions
$f$ on $P$ that change by $\langle u,\alpha\rangle$ under the action
of $\alpha\in\mathfrak{g}^{*}$ on $P$; the bracket on $\mathfrak{g}_{M}$ is
the Poison bracket on $P$.
\end{prop}
Here $\mathfrak{g}^{*}$ is a Poisson Lie group in the standard way (with
the Kirillov-Kostant Poisson structure $\pi_{\mathfrak{g}^{*}}$); a Poisson
principal $\mathfrak{g}^{*}$-bundle is a principal $\mathfrak{g}^{*}$-bundle
$P\rightarrow M$ with a Poisson structure on $P$ such that the action $P\times\mathfrak{g}^{*}\rightarrow P$
is a Poisson map.

Proof of the proposition is straightforward. If we choose a trivialization
of $P$, in terms of the decomposition \[
{\textstyle \bigwedge}^{2}T(M\times\mathfrak{g}^{*})={\textstyle \bigwedge}^{2}TM\oplus TM\otimes\mathfrak{g}^{*}\oplus{\textstyle \bigwedge}^{2}\mathfrak{g}^{*},\]
the Poisson structure on $P=M\times\mathfrak{g}^{*}$ is equal to\begin{equation}
\sigma^{0}+\sigma^{1}+(\sigma^{2}+\pi_{\mathfrak{g^{*}}}),\label{eq:decp}\end{equation}
where $\pi_{\mathfrak{g^{*}}}$ is the Poisson structure on $\mathfrak{g}^{*}$
and $\sigma=\sigma^{0}+\sigma^{1}+\sigma^{2}$ satisfies the MC equation
(\ref{eq:MC}), since $0=[\sigma+\pi_{\mathfrak{g^{*}}},\sigma+\pi_{\mathfrak{g^{*}}}]=2[\pi_{\mathfrak{g^{*}}},\sigma]+[\sigma,\sigma]=
2\,d\sigma+[\sigma,\sigma]$. Choosing a different trivialization corresponds to
a gauge transformation of $\sigma$.

As we noticed in Section \ref{sec:HPA}, a true Poisson $\mathfrak{g}$-action can be characterized as a 
$\mathfrak{g}$-HPA together with
a morphism to the trivial $\mathfrak{g}$-HPA on a point. In other words, a Poisson $\mathfrak{g}$-action is the
same as a principal Poisson $\mathfrak{g}^*$-bundle $P$ together with a $\mathfrak{g}^*$-equivariant Poisson map
$P\rightarrow\mathfrak{g}^*$.

Given a $\mathfrak{g}$-HPA on a Poisson manifold $M$ we cannot say directly what is the quotient of $M$ by the action.
There is, however, a well-defined ``quotient up to Morita equivalence'', namely the Poisson manifold $P$ (the total space
of the corresponding principal $\mathfrak{g}^*$-bundle).

Here is a somewhat exotic example of an up to homotopy Poisson action.
As was observed by A.~Vaintrob \cite{V}, if $A\rightarrow M$ is a vector bundle,
the structure of a Lie algebroid on $A$ is equivalent to a degree-1
vector field $Q$ on the graded supermanifold $A[1]$ such that $Q^{2}=0$,
and also to a degree $-1$ odd Poisson structure $\pi$ on $A^{*}[1]$;
a bialgebroid structure is then both $\pi$ and $Q$ on $A^{^{*}}[1]$
such that $\mathcal{L}_{Q}\pi=0$. Poisson actions up to homotopy
appear in quasi-bialgebroids: $A$ is a quasi-bialgebroid if we have
a degree $-1$ odd Poisson structure $\pi$ on $A^{*}[1]$ and an
up to homotopy action of the 1-dimensional graded Lie algebra generated
by $Q$ ($[Q,Q]=0$, $\deg Q=1$) on $A^{*}[1]$. In other words,
as was noticed in \cite{S1}, it is a Poisson principal $\mathbb{R}[2]$-bundle
$P\rightarrow A^{*}[1]$ (where $\mathbb{R}[2]$ has zero Poisson structure).

\section{Up to homotopy action of a group on an algebra\label{sec:HAA}}

An \emph{up to homotopy action} \emph{of a Lie algebra $\mathfrak{g}$ on
an associative algebra $A$} (a \emph{$\mathfrak{g}$-HAA on $A$} for short)
is an extension $\mathfrak{g}_{A}$ of $\mathfrak{g}$ by $A$ (where $A$ is
considered as a Lie algebra with the bracket given by the commutator),
such that for any $X\in\mathfrak{g}_{A}$, the map $\mathit{ad}_{X}:A\rightarrow A$
is a derivation of the associative algebra $A$. 

Similarly, an \emph{up to homotopy action} \emph{of a group $G$ on
an algebra} (a \emph{$G$-HAA}) is a $G$-graded algebra $A_{G}$,
i.e.~an algebra of the form $A_{G}=\bigoplus_{g\in G}A_{g}$ with
product mapping $A_{g}\otimes A_{h}$ to $A_{gh}$, such that every
$A_{g}$ contains an invertible element. The algebra $A_{e}$ is the
algebra on which $G$ `acts up to homotopy'. In this definition, $G$
is just an abstract group (and vector spaces are over an arbitrary
field). In a smooth definition, $G$ would be a Lie group and $A_{g}$'s
would form a smooth vector bundle over $G$.

\newcommand{\br}[1]{\langle#1 \rangle}

One can give a more combinatorial definition. Let us choose an invertible
element $\br{g}$ in every $A_{g}$. For every $g\in G$ we have an
automorphism $\rho(g)$ of $A_{e}$, given by $\rho(g)a=\br{g}a\br{g}^{-1}$.
This $\rho$ is not necessarily an action of $G$; instead, we have
\begin{equation}
\rho(g)\rho(h)a=c(g,h)(\rho(gh)a)c(g,h)^{-1},\label{eq:ha1}\end{equation}
where the elements $c(g,h)\in A_{e}$ are given by $c(g,h)=\br{g}\br{h}\br{gh}^{-1}$.
If we compute $\br{g_{1}}\br{g_{2}}\br{g_{3}}$ in two ways, we get
a cocycle identity \begin{equation}
c(g_{1},g_{2})c(g_{1}g_{2},g_{3})=(\rho(g_{1})c(g_{2},g_{3}))c(g_{1},g_{2}g_{3}).\label{eq:ha2}\end{equation}
Vice versa, when we are given an algebra $A_{e}$, its automorphisms
$\rho(g)$ and invertible elements $c(g,h)$ satisfying (\ref{eq:ha1})
and (\ref{eq:ha2}), we get an up to homotopy action.

This is indeed a non-commutative analogue of a principal Poisson $\mathfrak{g}^{*}$-bundle.
To get the analogue one replaces Poisson spaces by associative algebras
and reverses the arrows; the counterpart of $\mathfrak{g}^{*}$ is a group
algebra $\mathbb{F}[G]$ (where $\mathbb{F}$ is our base field).
The analogue of a Poisson space with a Poisson action of $\mathfrak{g}^{*}$
is then an algebra with a compatible structure of a $\mathbb{F}[G]$
comodule, i.e.~a $G$-graded algebra. The analogue of a principal
$\mathfrak{g}^{*}$-action is the existence of invertible elements in every
$A_{g}$; the algebra $A_{e}$ is the analogue of the base of the
principal bundle.

A true action of $G$ on an algebra $A_{e}$ should be, of course,
a special case of an up to homotopy action. All the spaces $A_{g}$
are in this case equal to $A_{e}$; to distinguish between them, let
us denote $A_{g}$ by $A_{e}\otimes g$ (its elements will be denoted
by $a\otimes g$, $a\in A_{e}$). The product $A_{g}\times A_{h}\rightarrow A_{gh}$
is given by $(a\otimes g)(b\otimes h)=(a\, g\cdot b)\otimes gh$.
This $G$-graded algebra is the crossed product of $A_{e}$ with $G$.

As an example of a different type, let $\tilde{G}$ be a central extension
of $G$ by $\mathbb{F}^{*}$, where $\mathbb{F}$ is our base field.
Then the line bundle associated to the principal $\mathbb{F}^{*}$-bundle
$\tilde{G}\rightarrow G$ is an up to homotopy action. And vice versa, every
up to homotopy action which is a line bundle (i.e.~where $A_{g}$'s
are 1-dimensional), is of this type; in other words, an up to homotopy
action of $G$ on $\mathbb{F}$ is equivalent to a central extension
of $G$ by $\mathbb{F}^{*}$.

\section{Poisson vs.~Hamiltonian actions: crossed products in Poisson geometry}

For any Poisson manifold $M$, a Poisson map $M\rightarrow\mathfrak{g}^{*}$ (a
moment map) gives rise to a Poisson action of $\mathfrak{g}$ on $M$; Poisson
actions of this form are called \emph{Hamiltonian}. For every Poisson
action one can construct a universal Hamiltonian action; more precisely,
the forgetful functor \begin{equation}
\textrm{Hamiltonian actions }\rightarrow\textrm{ Poisson actions}\label{eq:forg}\end{equation}
has a left adjoint. As we shall see, this adjoint is connected with
up to homotopy Poisson actions.

As we noticed in Section \ref{sec:HPA}, a Poisson action of $\mathfrak{g}$
on $M$ is the same as a $\mathfrak{g}$-HPA on $M$ and a morphism to the
trivial $\mathfrak{g}$-HPA on a point. Using the language of principal
Poisson bundles, a Poisson action of $\mathfrak{g}$ on $M$ is thus the
same as a principal Poisson $\mathfrak{g}^{*}$ bundle $P\rightarrow M$ with a
Poisson $\mathfrak{g}^{*}$-equivariant map $P\rightarrow\mathfrak{g}^{*}$. Since we
have a Poisson map $P\rightarrow\mathfrak{g}^{*}$, on $P$ we have a Hamiltonian
action of $\mathfrak{g}$. Thus for any Poisson action (on $M$) we can
find a Hamiltonian action (on $P=M\times\mathfrak{g}^{*}$); this operation
is left adjoint to the forgetful functor (\ref{eq:forg}), i.e.~it
satisfies a universal property:

\begin{prop}
The forgetful functor (\ref{eq:forg}) has a left adjoint given by
$M\mapsto P=M\times\mathfrak{g}^{*}$, with the Poisson structure given
by (\ref{eq:decp}) (with $\sigma^{0}$ the Poisson structure on $M$,
$\sigma^{1}$ the action of $\mathfrak{g}$ on $M$ and $\sigma^{2}=0$).
\end{prop}
Proof of this proposition is straightforward: if $N$ is a Poisson
manifold with a moment map $\mu:N\rightarrow\mathfrak{g}^{*}$, and if $f:N\rightarrow M$
is a $\mathfrak{g}$-equivariant Poisson map, we should provide a unique
Poisson map $\tilde{f}:N\rightarrow P$ commuting with the moment maps; the
map $\tilde{f}$ is simply $f\times\mu$.

Let us translate the previous to the world of noncommutative algebras.
There the analogue of a moment map is a product and unit preserving
map $\mu:G\rightarrow A$ (i.e. an algebra morphism $\mathbb{F}[G]\rightarrow A$);
it gives an action of $G$ on $A$ by $g\cdot a=\mu(g)a\mu(g^{-1})$.
The right-adjoint operation is the crossed product: given an action
of $G$ on $A$, it is $A\otimes\mathbb{F}[G]$ as a vector space,
with the product given by \[
(a_{1}\otimes g_{1})(a_{2}\otimes g_{2})=(a\, g_{1}\cdot a_{2})\otimes(g_{1}g_{2});\]
the map $\mu$ is given by $\mu(g)=1\otimes g$.

We can also generalize our proposition to so called non-equivariant
moment maps, i.e.~Poisson maps $M\rightarrow\mathfrak{p}$, where $\mathfrak{p}$ is
a principal Poisson $\mathfrak{g}^{*}$-homogeneous space. Such a map again
gives a Poisson action of $\mathfrak{g}$ on $M$. The adjoint functor is
the same as in the proposition: $M\mapsto M\times\mathfrak{p}$ (with the
mixed part of the Poisson structure equal to the action of $\mathfrak{g}$
on $M$). Notice that $\mathfrak{p}$ is a $\mathfrak{g}$-HPA on a point. A $\mathfrak{g}^{*}$
equivariant map $P\rightarrow\mathfrak{p}$ makes a principal Poisson $\mathfrak{g}^{*}$-bundle
$P\rightarrow M$ to a true Poisson action of $\mathfrak{g}$ on $M$.

The non-commutative analogue of this generalization uses a central
extension $\tilde{G}$ of a group $G$ by $\mathbb{F}^{*}$. A non-equivariant
moment map is a product preserving map $\mu:\tilde{G}\rightarrow A$ which
is the identity on $\mathbb{F}^{*}$. Such a map gives an action of
$G$ on $A$ via $g\cdot a=\mu(\tilde{g})a\mu(\tilde{g}^{-1})$, where
$\tilde{g}\in\tilde{G}$ is any element over $g\in G$. There is an
adjoint functor that is an obvious generalization of the crossed product.

\section{Quantization of up to homotopy Poisson actions}

The formality theorem of Kontsevich gives us a quantization of up
to homotopy Poisson actions depending on a formal parameter $\hbar$.
The resulting up to homotopy action
on the quantized algebra of functions will be also formally close
to a true action.

Let us thus have a $\mathfrak{g}$-HPA on $M$ formally depending on $\hbar$ and let
the corresponding Poisson structure on $M$ be $O(\hbar)$ (as is usual in deformation quantization).
For the lack of fantasy
we'll call such a HPA \emph{formally good}; these are the HPA's we're going to quantize.

When we decompose the HPA to $\sigma=\sigma^{0}+\sigma^{1}+\sigma^{2}$ as in (\ref{eq:dec}),
$\sigma^{0}$ is the Poisson structure on $M$ and hence is of the form
\[\sigma^{0}=\sigma_{(1)}^{0}\hbar+\sigma_{(2)}^{0}\hbar^{2}+\sigma_{(3)}^{0}\hbar^{3}+\dots;\]
$\sigma^1$ and $\sigma^2$ are not necessarily $O(\hbar)$, but from the MC equation (\ref{eq:MC})
we get that $\sigma_{(0)}^{1}$ is a true action of $\mathfrak{g}$ on $M$.

Since $\sigma=\sigma^{0}+\sigma^{1}+\sigma^{2}$ satisfies the MC
equation (\ref{eq:MC}), it can be directly plugged into the formality
theorem of Kontsevich, but only if $\sigma^{1}=O(\hbar)$. If we want
to allow $\sigma_{(0)}^{1}\neq0$, we need to use a rather straightforward
trick. Under our assumptions we'll get a $\mathfrak{g}$-HAA on the algebra
of quantized functions on $M$. If moreover the $\mathfrak{g}$-action $\sigma_{(0)}^{1}$
comes from an action of the simply-connected Lie group $G$ on $M$,
we'll also get a $G$-HAA. 

We shall use quantization of Poisson families from \cite{K} (and
terminology from \cite{S2}). Recall that an \emph{algebroid} over
a set $S$ is a linear category whose set of objects is $S$ and where
any two objects are isomorphic. In deformation quantization, algebroids
are produced by quantization of \emph{Hamiltonian families of formal
Poisson structures}. Such a family is by definition a solution of
the MC equation in the DGLA $\Omega(B)\hat{\otimes}L_{M}'$, where
$B$ and $M$ are manifolds, $\hat{\otimes}$ is the completed tensor
product (it just means that the coefficients are arbitrary smooth
functions on $B\times M$) and the DGLA $L_{M}'$ is given by\[
L_{M}'^{-1}=C^{\infty}(M)[[\hbar]]\]
\[
L_{M}'^{i}=\hbar\Gamma({\textstyle \bigwedge^{i+1}}TM)[[\hbar]].\]
If we choose a family of connections on $M$ parametrized by $B$,
we can quantize the Hamiltonian family to a {}``\emph{tight family
of $*$-products}'', i.e.~to a solution of the MC equation in $\Omega(B)\hat{\otimes}L_{M}''$,
where the DGLA $L_{M}''$ is given by\[
L_{M}''^{-1}=C^{\infty}(M)[[\hbar]]\]
\[
L_{M}''^{i}=\hbar\mathit{PD}^{i+1}[[\hbar]],\]
where $\mathit{PD}^{i+1}$ is the space of polydifferential operators\[
\underbrace{C^{\infty}(M)\times\dots\times C^{\infty}(M)}_{i+1\textrm{ times}}\rightarrow C^{\infty}(M).\]
As Kontsevich proved, this quantization is natural (i.e.~diffeomorphism-invariant)
and local. Moreover, if we suppose that $B$ is 2-connected, the quantized
family gives us naturally an algebroid over $B$.

Suppose now that $G$ is a group and that we are given a $G$-HAA
$A_{G}$. We then get an algebroid over $G$, defined by $\mathit{Hom}(g,h)=A_{g^{-1}h}$.
Vice versa, given an algebroid over $G$, in which we can identify
$\mathit{Hom}(g,h)$'s for a fixed $g^{-1}h$ in a way compatible
with the product, we get an up to homotopy action.

We shall now take for $G$ the 1-connected (and therefore also 2-connected)
Lie group integrating $\mathfrak{g}$, set $B=G$, construct a Hamiltonian
family, quantize it to an algebroid over $G$, and finally notice
that the family is right $G$-invariant, which will thus turn the
algebroid into a $G$-HAA.

The first thing is to notice the connection between Hamiltonian families
of formal Poisson structures and formally good $\mathfrak{g}$-HPAs:

\begin{lem}
Let $\rho$ be a (left) action of $\mathfrak{g}$ on $M$. A formally good
$\mathfrak{g}$-HPA $\sigma$ on $M$, such that $\sigma_{(0)}^{1}=\rho$,
is equivalent to a right-$\mathfrak{g}$-invariant Hamiltonian family of
formal Poisson structures on $M$ parametrized by $G$, where the
right $\mathfrak{g}$ action is by right translations on $G$ and by $-\rho$
on $M$.
\end{lem}
Suppose now that the $\mathfrak{g}$-action $\rho$ integrates to a $G$-action.
Then we can choose a right-$G$-invariant family of connections on
$M$ parametrized by $G$ (just choose one connection on $M$ and
then use the action). The algebroid over $G$ we get will be right-$G$-invariant,
so by the above remarks we get a $G$-HAA. We also get a $\mathfrak{g}$-HAA:
after quantization we have a right-$G$-invariant solution of the
MC equation in the DGLA $\Omega(G)\hat{\otimes}L_{M}''$.

If $\rho$ doesn't integrate to a $G$-action, we can choose the invariant
family of connections only locally. It is, of course, sufficient to
get a $\mathfrak{g}$-HAA. We also get a {}``local $G$-HAA'', though
we leave it to the reader to formulate the definition.

Let us summarize the outcome of this section:

\begin{prop}
Let $\rho$ be a (left) action of $\mathfrak{g}$ on $M$ and let $\sigma$
be a formally good $\mathfrak{g}$-HPA on $M$ with $\sigma_{(0)}^{1}=\rho$.
\begin{enumerate}
\item There is a quantization of $\sigma$ to a $\mathfrak{g}$-HAA on a quantized
algebra of functions on $M$
\item If moreover $\rho$ integrates to an action of $G$ (the 1-connected
group integrating $\mathfrak{g}$) then we also get a $G$-HAA.
\end{enumerate}
\end{prop}

\section{Generalizations for Lie bialgebras and Hopf algebras; Open ends}

There is now an obvious definition of up to homotopy actions of Lie
bialgebras. Let $(\mathfrak{g},\mathfrak{g}^{*})$ be a Lie bialgebra; let us
also choose a Poisson Lie group $G^{*}$ integrating $\mathfrak{g}^{*}$.
An \emph{up to homotopy Poisson action} of $(\mathfrak{g},G^{*})$ ($(\mathfrak{g},G^{*})$-HPA
for short) is, by definition, a principal Poisson $G^{*}$ bundle.
Such a bundle $P\rightarrow M$ with a $G^{*}$-equivariant Poisson map $P\rightarrow G^{*}$
is again equivalent to a Poisson action of $(\mathfrak{g},\mathfrak{g}^{*})$
on $M$: indeed, $P$ becomes $M\times G^{*}$, and the mixed term
of the Poisson structure is the action of $\mathfrak{g}$. Moreover,
this operation is left-adjoint to the forgetful functor\[
\textrm{Poisson maps to $G^{*}$ (moment maps) }\rightarrow\textrm{ Poisson actions of $(\mathfrak{g},\mathfrak{g}^{*})$}.\]
More generally, we have the following:

\begin{prop}
\label{pro:bialg-adj}Let $\mathfrak{P}$ be a right principal Poisson homogeneous
space of $G^{*}$. If $M$ is a Poisson manifold, a Poisson map $\mu:M\rightarrow\mathfrak{P}$
generates a (left) Poisson action $\rho$ of $(\mathfrak{g},\mathfrak{g}^{*})$
on $M$ via $\rho(v)_{x}=\langle\pi_{M},\mu^{*}v\rangle$, where $x\in M$,
$v\in\mathfrak{g}\cong T_{\mu(x)}^{*}\mathfrak{P}$ (these spaces are isomorphic
because $\mathfrak{P}$ is a principal $G^{*}$ space) and $\pi_{M}$ is
the Poisson structure on $M$.\footnote{%
Up to this point the proposition is of course well known and is included
just for completeness. }%
This functor has a left adjoint, given
by $M\mapsto P=M\times\mathfrak{P}$, with the Poisson structure $\pi_{P}=\pi_{M}+\pi_{\mathfrak{P}}+\textrm{the $\mathfrak{g}$-action}$.
Any Poisson principal $G^{*}$-bundle $P\rightarrow M$ with a Poisson equivariant
map $P\rightarrow\mathfrak{P}$ is of this form.
\end{prop}

Let $\mathfrak{d}$ be the Drinfeld double of $(\mathfrak{g},\mathfrak{g}^{*})$. Poisson
principal homogeneous $G^{*}$-spaces, or in other words Poisson principal
$G^{*}$-bundles over a point, are classified by Lagrangian subalgebras
of $\mathfrak{d}$ transversal to $\mathfrak{g}^{*}$ \cite{DS}. There is a similar
classification of general Poisson principal bundles. In the case of
the trivial bundle we have:

\begin{prop}
A Poisson structure on $P=M\times G^{*}$, making it to a Poisson
principal $G^{*}$ bundle, is the same as a Dirac structure in the
Courant algebroid $\mathfrak{d}\oplus(T\oplus T^{*})M$, transversal to
$\mathfrak{g}^{*}\oplus TM$.
\end{prop}
If we also require the projection $M\times G^{*}\rightarrow G^{*}$ to be
Poisson (i.e.~when we are interested in the true Poisson actions
of $(\mathfrak{g},\mathfrak{g}^{*})$ on $M$), the Dirac structure in $\mathfrak{d}\oplus(T\oplus T^{*})M$
has to project to $\mathfrak{g}\subset\mathfrak{d}$ under the projection $\mathfrak{d}\oplus(T\oplus T^{*})M\rightarrow\mathfrak{d}$.
This is a special case of an observation from \cite{BCS}: if $\mathfrak{g}\subset\mathfrak{d}$
is a Lie quasi-bialgebra, a $\mathfrak{g}\subset\mathfrak{d}$-quasi-Poisson structure
on a manifold $M$ (i.e.~a {}``quasi-Poisson action of $\mathfrak{g}$
on $M$'') is the same as a Dirac structure in $\mathfrak{d}\oplus(T\oplus T^{*})M$
that projects to $\mathfrak{g}\subset\mathfrak{d}$ under the projection $\mathfrak{d}\oplus(T\oplus T^{*})M\rightarrow\mathfrak{d}$,
and that is transversal to $TM$. (This suggests a definition of {}``up
to homotopy quasi-Poisson structures'', but it's not clear what it
would be good for.)

In the case of a general (non-trivial) principal $G^{*}$-bundle $P\rightarrow M$,
the Poisson structures on $P$ making it to a Poisson principal bundle,
are again the same as Dirac structures in some transitive Courant
algebroid $C\rightarrow M$, transversal to $A\subset C$, where $A$ is
the Atiyah Lie algebroid corresponding to $P$. The Courant algebroid
$C$ is obtained by first inducing $P$ to a principal $D$-bundle
(i.e.~by creating the associated principal $D$-bundle to $P$),
where $D$ is a Lie group integrating $\mathfrak{d}$, and then using
the reduction procedure from \cite{SL}. (Since the reduction procedure
is shown in \cite{SL} to produce basically all transitive Courant
algebroids, this again suggests some {}``up to homotopy and quasi''
directions). The fact that we get a Lie bialgebroid (a pair of transversal
Dirac structures in a Courant algebroid) is not surprising -- the
corresponding Poisson groupoid is the gauge groupoid of $P$.

Homotopy Poisson actions, as we defined them, involve a Lie algebra $\mathfrak{g}$. To integrate
a HPA to a kind of ``action'' of a Lie group $G$, we can take the Poisson principal bundle $P$,
form its gauge groupoid (that is automatically a Poisson groupoid), and finally pass
to its dual Poisson groupoid $\tilde \Gamma$ if it exists. The Poisson groupoid $\tilde\Gamma$ can be characterized
as follows: it is an extension of $G$ by a symplectic groupoid $\Gamma$ integrating the Poisson manifold $M$.
This kind of extension might be called ``up to homotopy action of $G$ on $M$''; notice that it requires $M$ to be integrable.
It is a direct translation of the definition of $G$-HAAs from Section \ref{sec:HAA}.

Let us now switch to the world of non-commutative algebras, replacing
Lie bialgebras by Hopf algebras. Let $H$ be a Hopf algebra (the case
of $H=\mathbb{F}[G]$ will lead back to $G$-HAA's). The substitute
for an action of a group $G$ up to homotopy on an algebra $A$ is
a $H$-Galois extension of $A$, i.e.~an algebra $B$, with a $H$-comodule
structure $c:B\rightarrow B\otimes H$, $b\mapsto b_{(1)}\otimes b_{(2)}$
that is also an algebra map, and an algebra map $i:A\rightarrow B$, such
that \[
i(A)=\{ b\in B|c(b)=b\otimes1\}\]
and such that the following map is an isomorphism:\[
B\otimes_{A}B\rightarrow B\otimes H,\quad x\otimes y\mapsto xy_{(0)}\otimes y_{(1)}.\]

I don't know if there is a theorem that would give a quantization
of a $(\mathfrak{g},G^{*})$-HPA on $M$ to a Hopf-Galois extension of a
quantized algebra of functions on $M$.

The non-commutative analogue of Proposition \ref{pro:bialg-adj} remains
valid. The analogue of $\mathfrak{P}$ is an $H$-Galois extensions of the
base field $\mathbb{F}$; in the case of $H=\mathbb{F}[G]$ this is
the same as a central extension of $G$ by $\mathbb{F}^{*}$.

\end{document}